\documentclass[12pt]{article}

\usepackage{mathptmx}
\usepackage{amsmath}
\usepackage{amsthm}
\usepackage{amssymb}
\def\pf{{\bf Proof. }}
\def\beq{\begin{equation}}
\def\eeq{\end{equation}}

\def\nd{\noindent}
\def\l{\lambda}

\def\g{\gamma}

\def\e{\epsilon}

\def\al{\alpha}
\def\bt{\beta}
\def\ro{\rho}

\def\t{\tilde}
\def\ra{\rightarrow}
\def\npc{of nonpositive sectional curvature}

\def\qrn{Ricci quasi-negatively curved}

\def\nd{\noindent}
\def\la{\langle}
\def\ra{\rangle}
\def\p{\partial}
\def\hf{\hfill{$\Box$}}
\def\<{\leq}
\def\>{\geq}

\newtheorem{thm}{Theorem}[section]
\newtheorem{lem}{Lemma}[section]

\newtheorem{rem}{Remark}[section]

\begin{document}
\title{\bf Compact Lie Group Actions on Closed Manifolds 
of Non-positive Curvature}
\author{{\sc Bin Xu}\\
{\small \sc Department of Mathematics, Johns Hopkins University, Baltimore MD 21218, USA}\\
{\small \it E-mail: bxu@math.jhu.edu}
}
\date{}
\maketitle
\begin{abstract}
\noindent A. Borel proved that, if a finite group $F$ acts effectively and continuously on a closed 
 aspherical manifold $M$ with centerless fundamental group $\pi_1(M)$, then a natural 
homomorphism $\psi$ from $F$ to the outer automorphism group ${\rm Out}\,\pi_1(M)$ of $\pi_1(M)$, 
called the associated abstract kernel, is a monomorphism. 
In this paper,  we investigate to what extent Borel's theorem holds 
for a compact Lie group $G$ acting effectively and smoothly on a particular orientable aspherical manifold $N$ 
admitting a Riemannian metric $g_0$ of non-positive curvature in case that $\pi_1(N)$ has a non-trivial center.
It turns out that if $G$ attains the maximal dimension equal to the rank of Center $\pi_1(N)$ and the metric $g_0$ is real analytic, then any element of $G$ defining a diffemorphism homotopic
to the identity of $N$ must be contained in the identity component $G^0$ of $G$.
Moreover, if the inner automorphism group of $\pi_1(N)$ is torsion free, then the associated abstract kernel 
$\psi: G/G^0\to {\rm Out}\,\pi_1(N)$ is a monomorphism. The same result holds for the non-orientable $N$'s under 
certain techical assumptions. Our result is an application of a theorem by 
Schoen-Yau (Topology, {\bf 18} (1979), 361-380) on harmonic mappings. 
 
\end{abstract}

\nd {\small \bf Mathematics Subject Classifications (2000):} Primary: 57S25, 53C43;  Secondary: 20F34.  \\

\nd{\small \bf Key words:} aspherical manifold, compact Lie group action, non-positive curvature, harmonic mapping.

\section{Introduction}
\quad Let $M$ be a closed aspherical manifold. A lot of information has been obtained by 
Conner-Raymond \cite{CR1}-\cite{CR3} regarding group actions on 
such manifolds. One of their main results in \cite{CR1} is that a connected compact Lie group acting effectively and continuously
on $M$ is a torus group of dimension not greater than the rank of Center $\pi_1(M)$
(note that $\pi_1(M)$ is torsion free by the Smith fixed point theorem). In particular, any compact Lie group acting effectively and continuously
on a closed aspherical manifold with centerless fundamental group must be a finite group. 
A group action on $M$ is {\it homotopically trivial} if 
any element of this group, which defines a homeomorphism of $M$, 
is homotopic to the identity. Gottlieb-Lee-{\" O}zaydin 
\cite{GL} showed that a compact Lie group acting effectively
and homotopically trivially on $M$ is abelian. 
A simple but nontrivial example of 
Gottlieb-Lee-{\" O}zaydin's
result, is  the group $\la\exp ({2\pi \sqrt{-1}/m})\ra
\cong {\bf Z}/m{\bf Z}$  acting on
the circle $\displaystyle
{S^1=\{\exp (\sqrt{-1}\theta):\theta\in [0,\ 2\pi)\}}$ by multiplication.

Let $G$ be a compact Lie group acting effectively and continuously on  $M$ and $G^0$ be the 
identity component of $G$. We quickly review the definition of the {\it associated abstract kernel} that is
 a homomorphism $\psi:G\to {\rm Out}\,\pi_1(M)$ (cf Sections 2-3 \cite{CR3} for the detail). Choose any base point
$x$ in $M$ and for each $\al$ in $G$ select a path ${\mathcal P}_\al(t)$ joining $x$ to $\al x$. 
Define the corresponding automorphism $\al_*$ on $\pi_1(M,x)$ as follows: if $\sigma(t)$ is a loop
based at $x$, then $\al_*(\sigma)$ is represented by
\[\left\{
\begin{array}{ll}
{\mathcal P}_\al(3t) & 0\leq t\leq 1/3; \\
\al\bigl(\sigma(3t-1)\bigr) & 1/3\leq t\leq 2/3;\\
{\mathcal P}_\al(3-3t) & 2/3\leq t\leq 1.
\end{array}
\right.
\]
This $\al_*$ is unique up to an inner automorphism of $\pi_1(M,x)$ and yields a homomorphism 
$\psi:G\to {\rm Out}\,\pi_1(M)$, called the associated abstract kernel. Since $G^0$ is contained
in the kernel of this homomorphism, $\psi$ can be reduced to another homomorphism from $G/G^0$ to
${\rm Out}\,\pi_1(M)$, which we still call the associated abstract kernel and denote by $\psi$. 
A. Borel proved that, if $G$ is finite and $M$ has centerless fundamental group, then 
the associated abstract kernel $\psi:G\to {\rm Out}\,\pi_1(M)$ is a monomorphism
(cf Theorem 3.2 \cite{CR3}). In particular,
the two homeomorphisms defined by any two elements of $G$ are not homotopic.

Lawson and Yau \cite{LY} showed that the isometry group $I(N)$
of a closed Riemannian manifold $N$ of non-positive curvature  has dimension equal to 
the rank of Center $\pi_1(N)$, and that the identity component $I^0(N)$ of $I(N)$
is a torus group which is generated by the parallel vector fields on $N$ 
(in other words, every Killing vector field on $N$ is parallel). Obviously $I^0(N)$ acts 
homotopically trivially on $N$. But Lawson and Yau \cite{LY} 
proved that any homotopically trival element of $I(N)$ must be contained in $I^0(N)$.

Let $N$ be a closed smooth manifold admitting a Riemannian metric $g_0$ of non-positive
curvature and $G$ be a compact Lie group acting smoothly and effectively on $N$. If $\pi_1(N)$ is
centerless, then it follows from the results by Conner-Raymond and Borel that $G$ must be finite and the
associated abstract kernel $\psi:G\to {\rm Out}\,\pi_1(N)$ is a monomorphism. In this paper, 
we shall investigate the following two problems:\\

\nd {\it Problem 1.}\quad When should any homotopically trivial element of $G$ be contained in $G^0$ 
(we are looking for a stronger version of Gottlieb-Lee-{\" O}zaydin's result for our case)? \\

\nd {\it Problem 2.}\quad To what extent Borel's theorem holds for 
the associated abstract kernel $\psi:G/G^0\to {\rm Out}\,\pi_1(N)$ in case that 
$\pi_1(N)$ has a non-trivial center? \\

\nd It turns out that these two problems are interconnected. By the result of Lawson-Yau,
we are not missing any positive results if we assume that $G$ attains the maxiaml dimension equal to rank of 
Center $\pi_1(N)$. Since the fact that the inner automorphism group of $\pi(M)$ is torsion free is critical 
in the proof of Borel's theorem (cf \cite{CR3}), it is naturally for us to assume this property for $\pi_1(N)$ 
while considering Problem 2.

\begin{thm}
\label{thm:np}
Let $N$ be a closed smooth manifold admitting a real analytic Riemannian metric $g_0$ of non-positive
curvature and $G$ be a compact Lie group acting smoothly and effectively on $N$
and having positive dimension equal to the rank of {\rm Center}\,$\pi_1(N)$. 
Suppose that there exists no free smooth action of odd prime order on $N$ if 
$N$ is non-orientable.
Then any homotopically trivial element of $G$ must be contained in $G^0$. Moreover, if 
the inner automorphism group of $\pi_1(N)$ is torsion free, then the associated abstract kernel 
$\psi:G\to {\rm Out}\,\pi_1(N)$ is a monomorphism. 
\end{thm}

\begin{rem}
\label{rem:kl}
{\rm The previous theorem only gives an incomplete answer to those two problems for non-orientable $N$'s. 
For an example, 
the only non-orientable closed 2-manifold whose fundamental group has non-trivial center is the Klein bottle $K$.
Since there exists no non-trivial free action on $K$, we know the answer of Problem 1 for the Klein bottle $K$ is Yes. 
However, the author does not know the answer for $K\times S^1$. 
Since there exists an involutive inner automorphism of 
$\pi_1(K)$,  the theorem can not give the answer of Problem 2 for the Klein bottle $K$, either. 
But the author conjectures that the assumption for non-orientable $N$'s is removable.
}
\end{rem}

If a compact Riemannian manifold has non-positive  curvature and has a point at which
all sectional curvatures are negative, then  Sampson in Theorem 5 of \cite{Sam} used his rigidity theorem 
(cf Theorem 4, \cite{Sam}) of harmonic mappings to show that the isometry group of the Riemannian manifold
is finite, and no two of its elements are homotopic. Here we call a harmonic map is {\it rigid} if it is the
only harmonic map in its homotopy class. Due to Frankel \cite{Fr}, the same result as Sampson's also 
holds for a compact Riemannian manifold of non-positive  curvature and of negative definite Ricci 
tensor. We shall not only unify the results of Sampson and Frankel from the point of view of
Borel's theorem, but also give a stronger result than theirs in the following 

\begin{thm}
\label{thm:qn}
Let $(V, g_1)$ be a closed  Riemannian manifold of non-positive  curvature and having a point 
at which the Ricci tensor is negative definite. Then the following two statements hold{\rm :}\\
{\rm a.}  $\pi_1(V)$ is centerless.\\
{\rm b.}  Any compact Lie group $G$ acting effectively and smoothly on $V$ is finite and 
the associated abstract kernel $\psi:G\to {\rm Out}\,\pi_1(N)$ is a monomorphism.\\
\end{thm}

\nd{\sc A quick proof of Theorem \ref{thm:qn}:} 
It follows from a theorem of Bochner \cite{B} (also see pages 324-326 in \cite{Wu}) that every Killing field on $V$ is zero since
$V$ has non-positive Ricci curvature and has a point at which the Ricci tensor is negative definite. 
Hence the isometry group of $(V,g_1)$ is finite. 
By the result of Lawson and Yau, $\pi_1(N)$ is centerless. Then by the result of Conner-Raymond, the compact Lie group $G$
must be a finite group even if $G$ only acts continuously on $N$. By Borel's theorem the associated abstract kernel 
$\psi:G\to {\rm Out}\,\pi_1(N)$ is a monomorphism.
\\

This paper is organized as follows. In Section 2 
we cite the rigidity theorem of Schoen-Yau \cite{SY} and prove a rigidity theorem 
(cf Theorem \ref{thm:rig}) of 
the harmonic mappings to compact \qrn\ manifolds. 
The first of them will be used to prove Theorem \ref{thm:np},
and the second will be applied to an analytic proof of Theorem \ref{thm:qn}.
In Section 3, we prepare two lemmata for the proof of these two theorems. 
The last section of this paper consists of the proof of Theorem \ref{thm:np} and 
the analytic proof of Theorem \ref{thm:qn}.

\section{Rigidity theorems of harmonic mappings}
\nd {\bf Theorem SY} 
(cf Theorem 4 in Schoen-Yau \cite{SY})
{\it Suppose $M,\ N$ are compact connected real analytic 
Riemannian manifolds and $N$ has nonpositive 
 curvatures.
Suppose $h:M\to N $ is a surjective harmonic map and its induced map $h_*:\pi_1(M)\to \pi_1(N)$ 
is also surjective. Then the space of surjective harmonic maps homotopic to $h$ is represented 
by $\{\bt\circ h| \bt\in I^0(N)\}$.}

\begin{thm}
\label{thm:rig}
Let $h_0:M\to N$ be a harmonic mappings, where $M$ is compact
and $N$ \npc\ . Suppose that there exists a point $p$ in $M$ 
such that the followings hold:\\
{\rm (a)} at $h_0(p)$ the Ricci tensor of $N$ is negative definite,\\
{\rm (b)} the differential map $dh_0(p)$ of $h_0$ at $p$ is surjective.\\
Then $h_0$ is the only harmonic mapping in its homotopy class.
\end{thm}

Before the proof of the above theorem, we make a quick review on 
the formula for second variation of the energy (cf J. Jost \cite{Jo})
for a family of harmonic mappings.

Let $M$ be a compact, and $N$ a complete 
Riemannian manfiolds of dimension $m$ and $n$, 
respectively. In local coordinates, the metric tensor of $M$ is
written as $(\g_{\al\bt})_{\al,\bt=1,\cdots,m}$, 
and the one of $N$ as $(g_{ij})_{i,j=1,\cdots,n}$. 
We shall use the notation 
$(\g^{\al\bt})_{\al,\bt=1,\cdots,m}=
(\g_{\al\bt})^{-1}_{\al,\bt=1,\cdots,m}$ (inverse metric tensor).
Let $f:M\to N$ be a smooth map and  $f^{-1}TN$ the pullback bundle on $M$ by 
$f$ of the tangent bundle $TN$ of $N$. $f^{-1}TN$ has the metric 
$(g_{ij}(f(x))$, and the cotangent bundle $T^*M$ of course has 
the metric $(\g^{\al\bt})$.
Then the energy density $e(f)$ is defined as 
$\frac{1}{2}||df||^2$, 
which is the square of the norm of the differential of $f$ as a section
of the Riemannian bundle $T^*M \otimes f^{-1}TN$. The {\it energy} 
of $f:M\to N$ is $E(f):=\int_M e(f) dM$,
with $dM$ the volume form of $M$. The smooth map $f:M\to N$ is {\it harmonic}
if and only if it is a critical point of the energy functional $E$.

Let $F:M\times (-\e, \e)\to N,\ F_t(x)=F(x,t)$
be a family of smooth maps between Riemannian manfiolds,
in which $F_0(x)=F(x,0)=f$.
Then $\displaystyle{W:=\frac{\partial F}{\partial t}|_{t=0}}$ is a section of
$f^{-1}TN$. Let $\nabla$ denote the Levi-Civita connection in $f^{-1}TN$ and
$R^N$ the curvature tensor of $N$. Then 
we have the following\\

\nd {\bf Fact E} For the second
variation of energy the equality
\[\frac{\partial^2}{\partial t^2} E(F_t)|_{t=0}
=\int_M ||\nabla W||^2_{f^{-1}TN}\,dM
-\int_M {\rm trace}_M \la R^N(df,W)W,df\ra_{f^{-1}TN}\,dM \]
holds provided that $F(x,\cdot)$ is geodesic for every $x$.\\

Then we recall a result on homotopic harmonic mappings by Hartman \cite{Ha}.
\\

\nd {\bf Fact H} (cf \cite{Jo})  Assume that $N$ is a complete manifold \npc.\ 
Let $f_0, f_1:M\to N$ be homotopic harmonic mappings. Then
there exists a family $f_t:M\to N,\ t\in [0,\ 1]$, of harmonic mappings
connecting them, for which the energy $E(f_t)$ is independent of $t$, and 
for which every curve $\g_x(t):=f_t(x)$ is geodesic, and 
$||\frac{\p}{\p t}\g_x(t)||$ is independent of $x$ and $t$.
\\

\nd {\sc Proof of Theorem \ref{thm:rig}}\ \  Let $h_1:M\to N$ be a harmonic map
homotopic to $h_0$. We can find a family $h_t:M\to N,\ t\in [0,\ 1]$, 
of harmonic mappings connecting them with the property as Fact H. By Fact E,
since $N$ has nonpositive  curvature,
\begin{eqnarray*}
0&=&\frac{d^2}{dt^2} E(h_t)\\
&=&\int_M\biggl(||\nabla \frac{\p}{\p t} \g_x(t)||^2-
{\rm trace}_M \la 
R^N\bigl(dh_t(x),\frac{\p}{\p t} \g_x(t)\bigr)\frac{\p}{\p t} \g_x(t),dh_t(x)\ra_{h_t^{-1}TN}\biggr)
dM\\
&\geq& 0 \ .
\end{eqnarray*}
Hence we obtain that for any $x\in M$ and any $t\in [0,\ 1]$,
\[{\rm trace}_M \la 
R^N\bigl(dh_t(x),\,\frac{\p}{\p t} \g_x(t)\bigr)\frac{\p}{\p t} \g_x(t),\,dh_t(x)\ra_{h_t^{-1}TN}=0,\]
in particular,
\[{\rm trace}_M \la 
R^N\bigl(dh_0(p),\frac{\p}{\p t} \g_p(t)\bigr)
\frac{\p}{\p t} \g_p(t),dh_0(p)\ra_{h_0^{-1}TN}
=0\ .\]
By assumptions (a), (b) and the following Fact T, we can see that
the tangent vector $\frac{\p}{\p t} \g_p(t)$ of the geodesic $\g_p(t)$ vanishes. 
By Fact H every geodesic $\g_x(t),\ x\in M,\ t\in [0,\ 1]$
degenerates into a point. That is, $h_0=h_1$. \hf \\

\nd {\bf Fact T} Let $R$ be a $n\times n$ real symmetric matrix having
 nonpostive eigenvalues and negative trace. Suppose $B$ is a $m\times n$
 real matrix with rank $n$. Then the matrix $BRB^t$ also has negative 
 trace, where  $B^t$ is the transposed matrix of $B$. \\
 
\nd {\sc Proof} We may assume 
$R$ to be the diagonal matrix ${\rm diag}(\l_1,\cdots,\l_n)$, 
in which $\l_1<0$ and $\l_2,\cdots,\l_n\leq 0$, without loss of generality.
We denote $B=(b_{ij})_{1\leq i\leq m,\ 1\leq j\leq n}$. Since
$B$ has rank $n$, there exists a nonzero element $b_{i1}$ in the first row of 
$B$. By computation, we know 
${\rm trace}\ BRB^t\leq \l_1 b_{i1}^2<0$.\hf

\section{Two lemmata}
\begin{lem}
\label{lem:deg}
Let $M$ and $N$ be compact connected smooth Riemannian manifolds of the same dimension and $f:M\to N$
a smooth mapping. Assume deg\ $f=m\not= 0$ if $M$ is orientable, and
deg\ $f\equiv 1\pmod{2}$ if $M$ is non-orientable. Define a subgroup $A$
of $I(M)$ by 
\[A=\{ \al\in I(M)|f\circ\al=f\}.\]
Then if $M$ is orientable, the order of $A$ divides $m$. In particular,
in case of deg\ $f=\pm 1$, the group $A$ is trivial. If $M$ is non-orientable,
the order of $A$ is an odd integer.  
\end{lem}
\nd \pf
We firstly prove the part in which $M$ is orientable. Taking a regualr value
$y_0\in N$ whose preimages under $f$ are $x_1,\cdots,x_{2k+m}$, then
\[m=\mbox{deg}\ f=\Sigma_{i=1}^{2k+m} \mbox{sgn\ det}\ Jf(x_i),\]
where $Jf(x_i)$ is the Jacobian matrix of $f$ at point $x_i$. Without loss of 
generality, we set $m\geq 1$ and
\[\mbox{sgn\ det}\ Jf(x_i)=1,\ \mbox{for}\ 1\leq i\leq k+m;\]
\[  \mbox{sgn\ det}\ Jf(x_j)=-1,\ \mbox{for}\ k+m+1\leq j\leq 2k+m.\]
The group $A$ acts on the set $f^{-1}(y_0)$ by the definition of $A$.
We claim that $A$ acts freely on this set. In fact if an element
$\bt\in A$ has $x_1$ as its fixed point, then its differential at point $x_1$
is the identity map since $df(x_1)\circ d\bt(x_1)=df(x_1)$. 
Since $M$ is connected, $\bt$ must be the identity of $M$.

We also claim that $A$ should preserve the orientation of $M$. Otherwise, letting ${\bar A}$
be the subgroup of $A$ whose elements preserve the orientation of $M$, we have
$[A:{\bar A}]=2$. If we set
${\bar A}=\{g_1,\cdots,g_n\},\ A=\{g_1,\cdots,g_n,gg_1,\cdots,gg_n\}$, 
then the equalities
\[\mbox{sgn\ det}\ Jf(x)=\mbox{sgn\ det}\ Jf(g_i(x)),\ \mbox{sgn\ det}\ Jf(x)=-\mbox{sgn\ det}\ Jf(gg_i(x))\]
hold for any $x\in f^{-1}(y_0)$ and $1\leq i\leq n$. As $A$ acts freely on $f^{-1}(y_0)$, we obtain deg\ $f$=0
and a contradiction. Since $A$ preserves the orientation of $M$ and acts freely on the set 
$\{x_1,\cdots,x_{2k+m}\}$, $A$ should act freely on the two set $\{x_1,\cdots,x_{k+m}\}$ and 
$\{x_{k+m+1},\cdots,x_{2k+m}\}$ respectively so that the order of $A$ divides $m$.

When $M$ is non-orientable, we know the number of the set $f^{-1}(y_0)$ is odd because of 
deg\ $f\equiv 1 \pmod{2}$. By the same way, we can show that $A$ acts freely on $f^{-1}(y_0)$
so that the order of $A$ is also odd.\hf

\begin{rem}
{\rm Schoen and Yau proved a result (cf. Theroem 4 (i) of \cite{SY}) 
more general than Lemma \ref{lem:deg}. 
However we prefer a simple proof as above for our special case.}
\end{rem}

We also need a well known result as follows: 

\begin{lem}
\label{lem:lift}
Let $M$ be a non-orientable smooth manifold and $M'$ its orientable double
covering. Then 
a diffeomorphism of $M$ can be lifted to that of $M'$.
\end{lem}

\section{Proof of Theorems 1.1 and 1.2}
 
 \nd {\sc Proof of Theorem \ref{thm:np}} Let $r>0$ be the rank of Center $\pi_1(N)$. 
Then $I^0(N,\,g_0)$ is a torus group of dimension $r$.  
 Since $G$ is a compact Lie group acting smoothly on $N$, we can choose
 a real analytic Riemannian metric $g$ such that $G$ acts on $(N, g)$ isometrically.
 By Eells-Sampson \cite{ES}  We have a harmonic mapping $h:(N,g)\to (N,g_0)$ homotopic to the
identity. Let $H$ be the subgroup of $G$ consisting of elements homotopic to the identity map of $N$.
Since $H$ contains $G^0$, $H$ is a closed Lie subgroup of $G$ having the same dimension as $G$.
For an element $\al$ of $H$,  since $h$ and $h\circ\al$ are two homotopic harmonic mappings satisfying
 the condition of Theorem SY, there exists a unique element $\bt$
 of $I^0(N,g_1)$ such that  $h\circ\al=\bt\circ h$, which leads to a Lie group homomorphism 
 \[\ro:H\to I^0(N,g_1),\ \ro(\al):=\bt\ .\] 
We shall prove $H=G^0$ in the following.\\

\nd {\it Case 1}\quad If $N$ is orientable, then we have a monomorphism 
$\ro: H\to I^0(N,g_0)$ by Lemma \ref{lem:deg}. Remember that $I^0(N, g_0)$ is a torus group of
 dimension equal to rank of Center $\pi_1(N)$. 
 Sine $\dim G=\dim\ I^0(N,g_1)$, $\ro$ must be an isomorhism so that
$H$ is a torus group. But $H$ contains the tours group $G^0$ with the same dimension, 
so $H=G^0$. \\

\nd{\it Case 2}\quad Let $N$ be non-orientable. Suppose that Kernel $\ro$ is non-trivial. 
By Lemma \ref{lem:deg} there exists an element $\al$ of Ker $\ro$ whose order is an odd prime number. 
If $\al$ has a fixed point $x$ on $N$, by Corollary 6.2 in \cite{CR1} $\al_*$ induces a non-trivial
automorphism of $\pi_1(N,x)$, which contradicts that $\al$ is homotopic to the identity map of $N$.
Therefore, $\al$ acts freely on $N$, which contradict the assumption that there exists no free smooth action
of odd prime order on $N$.\\ 

By now we have showed that any homotopically trivial element of $G$ is contained in $G^0$.
Finally we prove the associate abstract kernel $\psi:G/G^0\to {\rm Out}\,\pi_1(N)$ is a monomorphism
if the inner automorphism group ${\rm Inn}\,\pi_1(N)$ is torsion free. 
By previous argument, we need only to show that an element $[\al]$ in $G/G^0$ inducing an inner automorphism 
of $\pi_1(N)$ is homotopically trivial. Indeed, since $[\al]$ is of finite order and the inner automophism group
of $G$ is torsion free, $\psi([\al])$ must be the identity map on $\pi_1(N)$. Since $N$ is a $K(\pi_1(N),1)$-space, 
$\al$ is homotopic to the identity map of $N$.\hf\\ 

\nd
{\sc Proof of Theorem \ref{thm:qn}} At first we prove the finiteness
of $G$. We only need to show that the identity component $G^0$ is
trivial.  Since $G$ is a compact Lie group acting on $N$, we can choose
 a Riemannian metric $g$ such that $G$ action on $(V,g)$ is isometric.
 By Eells-Sampson \cite{ES}, there exists a harmonic mapping
 $h:(V,g)\to (V, g_1)$ homotopic to the identity so that
 $h$ is surjective. We shall prove\\
 
 \nd {\bf Claim 1} {\it $h$ is the only harmonic mapping in its 
 homotopy class.}
 
 \nd {\small \sc Proof of Claim 1}  Since $(V, g_1)$ has a point at which the Ricci tensor is negative definite, 
there exists an open subset $U$ of $N$  on which  $(V, g_1)$ has negative definite Ricci tensor. 
By Sard's theorem and the surjectivity of $h$, we know the  regular value set of $h$ is a dense subset of $N$ 
so that  there exists a regular point of $h$ which is mapped into $U$.
Then we apply Theorem \ref{thm:rig}.\\
 
 By Claim 1 $G^0$ is contained in the set
 $\{\al\in G: h\circ\al=h\}$,  which is a finite set by Lemma \ref{lem:deg}. That is, $G^0$ is trivial.
 
 Then we prove that no two elements of $G$ are homotopic.
 We only need to prove that an element $\al$ of $G$ homotopic to the 
 identity is the identity.
 If $V$ is orientable, then by Claim 1 and Lemma \ref{lem:deg}
 it follows from the fact that $h$ has degree 1. In the following we may assume
 that $V$ is non-orientable. Let $V'$ be its orientable double covering
 space and $\g:V'\to V'$ the deck transformation. By Lemma \ref{lem:lift}
 $\al$ has two liftings ${\t \al},\g\circ{\t \al}$ on $V'$, in which
 ${\t \al}$ is homotopic to ${\rm id}_{V'}$. Since the order of ${\t \al}$
 is equal to that of $\al$ $\leq |G|$, applying the previous argument to the finite 
 group action on $V'$ generated by ${\t \al}$,  we can see that ${\t \al}$ is ${\rm id}_{V'}$ and then
 $\al$ also should be the identity of $V$.\\

\nd{\bf Claim 2} {\it $\pi_1(V)$ is centerless.}

\nd{\small \sc Proof of Claim 2} It has been proved in Introduction. We give an alternative proof as follows.
The previous argument tells us that $I(V,g_1)$ is a finite group. 
By the result of Lawson and Yau mentioned in Introduction, the center of $\pi_1(V)$ has rank zero. 
The proof is completed by the well known fact that $\pi_1(V)$ is torsion free. \\

Finally we show that the abstract kernel $\psi:G\to {\rm Out}\,\pi_1(V)$ is a monomorphism.  
We have to prove that an element $\al$ of $G$ inducing an inner automorphism $\iota(\al)$
of $\pi_1(V)$ must be the identity map of $V$. By Claim 2, the inner automorphism group ${\rm Inn}\,\pi_1(V)$ 
of $\pi_1(V)$ is isomorphic to $\pi_1(V)/{\rm Center}\,\pi_1(V)=\pi_1(V)$, which is torsion free. Since $\al$ 
is of finite order and ${\rm Inn}\,\pi_1(V)$ is torsion free, the inner automorphism $\iota(\al)$ 
induced by $\al$ must be the identity map of $\pi_1(V)$. 
Since $V$ is a $K(\pi_1(V),1)$-space, the diffeomorphism $\al$ on $V$ must be homotopic to the identity map of $V$.  
The previous argument gives us that $\al$ should equal the identity map of $N$.\hf\\

\nd
{\bf Acknowledgements}

\nd This study was supported in part by the Japanese Government Scholarship, the JSPS Postdoctoral
Fellowship for Foreign Researchers and the FRG Postdoctoral Fellowship.


\begin{thebibliography}{99}
\bibitem{B}{\sc S. Bochner,} 
{\it Vector Fields and Ricci Curvature},
Bull. Amer.  Math. Soc. {\bf 52} (1946), 776-797.
 
\bibitem{CR1}{\sc P.\ Conner and F.\ Raymond,}
Actions of Compact Lie Groups on Asperical Manifolds,
in {\it Topology of Manifold} (Proceedings of the University of Georgia Topology of Manifolds Institute), 
ed. J. C. Cantrell and C. H. Edwards, Jr.  Markham, (1970), 227-264

\bibitem{CR2}{\sc P.\ Conner and F.\ Raymond,}
Injective Operations of the Toral Groups, {\it Topology}, {\bf 10} (1971),
283-296

\bibitem{CR3}{\sc P.\ Conner and F.\ Raymond,}
Manifolds with Few Periodic Homeomorphisms,
in {\it Lecture Notes in Math.} {\bf 299}, ed. H. T. Ku, L. N. Mann, J. L. Sicks and J. C. Su, 
Springer-Velag, New York, (1972),  1-75

\bibitem{ES}{\sc James Eells Jr and J.\ H.\ Sampson,}\ 
Harmonic Mappings of Riemannian Manifolds, {\it Amer. J. Math.}  
{\bf 86} (1964), 109-160

\bibitem{Fr}{\sc T. T. Frankel,}
On Theorem of Hurwitz and Bochner,
{\it  J. Math. Mech.} {\bf 15} (1966), 373-377 


\bibitem{Ha}{\sc P. Hartman},
On Homotopic Harmonic maps,
{\it Can. J. Math.} {\bf 19} (1967), 673-687


\bibitem{GL}{\sc D.\ H.\ Gottlieb and K.\ B.\ Lee and M.\ {\" O}zaydin,}\ 
Compact Group Actions and Maps into $K(\pi,1)$-Space,\ 
{\it Trans.  Amer. Math. Soc}  {\bf 287} (1985), 419-429

\bibitem{Jo}{\sc J. Jost,}
{\it Riemannian Geometry and Geometric Analysis},
3rd edition, Springer, 2000.


\bibitem{LY}{\sc H.\ B.\ Lawson and S.\ T.\ Yau,}\  Compact Manifolds with 
Nonpositive Curvature, {\it J. Diff. Geo}, {\bf 7}
\ (1972),  211-228

\bibitem{Sam}{\sc J.\ H.\ Sampson,}\ Some Properties and Applications of 
Harmonic Mappings, {\it  Ann.\ {\' E}cole Norm.\ Sup.} (4) {\bf 11} (1978), no. 2,
211-228

\bibitem{SY}{\sc R.\ Schoen and S.\ T.\ Yau,}\ Compact Group Actions 
and the Topology of Manifolds with Nonpositive Curvature,
{\it Topology}, {\bf 18} (1979), 361-380

\bibitem{Wu}{\sc Hung-Hsi Wu}, 
The Bochner Technique in Differential Geometry, 
{\it Math. Rep.} {\bf 3} (1988), no. 2, 289-3538. 
\end{thebibliography}
\end{document}